\documentclass[10pt]{amsart}
\usepackage{amsmath}
\usepackage{cases}
\usepackage{mathrsfs}
\usepackage{bbm}
\usepackage{amssymb}
\usepackage{amscd}
\usepackage{amsfonts,latexsym,amsmath,amsthm,amsxtra,mathdots,amssymb,latexsym,mathabx}
\usepackage[all,cmtip]{xy}
\RequirePackage{amsmath} \RequirePackage{amssymb}
\usepackage{color}
\usepackage{colordvi}
\usepackage{multicol}
\usepackage{hyperref}
\usepackage{mathtools}
\usepackage[margin=1.1in]{geometry}
\usepackage{xcolor}
\hypersetup{
    colorlinks,
    linkcolor={red!50!black},
    citecolor={blue!50!black},
    urlcolor={blue!80!black}
}

\newtheorem{Theorem}{Theorem}[section]
\newtheorem{Lemma}[Theorem]{Lemma}
\newtheorem{Corollary}[Theorem]{Corollary}

\newtheorem{Remark}[Theorem]{Remark}
\newcommand{\sumstar}{\sideset{}{^\star}\sum}

\begin{document}

\title{Triple correlations of Fourier coefficients of cusp forms}

\author[Y. Lin]{Yongxiao Lin}
\address{Department of Mathematics, The Ohio State University\\ 231 W 18th Avenue\\
Columbus, Ohio 43210-1174}
\email{lin.1765@osu.edu}

\begin{abstract}
We treat an unbalanced shifted convolution sum of Fourier coefficients of cusp forms. As a consequence, we obtain an upper bound for correlation of three Hecke eigenvalues of holomorphic cusp forms $\sum_{H\leq h\leq 2H}W\big(\frac{h}{H}\big)\sum_{X\leq n\leq 2X}\lambda_{1}(n-h)\lambda_{2}(n)\lambda_{3}(n+h)$, which is nontrivial provided that $H\geq X^{2/3+\varepsilon}$. The result can be viewed as a cuspidal analogue of a recent result of Blomer \cite{blo} on triple correlations of divisor functions.
\end{abstract}



\maketitle

\section{Introduction}

Recently Blomer \cite{blo} established an asymptotic formula with power saving error term for certain types of triple correlations of divisor functions. Motivated by his work, we are going to prove a cuspidal analogue. While shifted convolutions of two Fourier coefficients have been an extensively studied topic, there seem to be few results available on the correlation of three Fourier coefficients, with power saving error term. For instance, one of the highly interesting and open problems in analytic number theory is to find an asymptotic formula for 
\begin{equation}
D_h(X)=\sum_{n\leq X}\tau(n-h)\tau(n)\tau(n+h),
\end{equation}
where $\tau(n)$ denotes the divisor function. It is conjectured that 
\begin{equation}\label{D_h}
D_h(X)\sim c_h X(\log X)^3,
\end{equation}
as $X\rightarrow \infty$, for some positive constant $c_h$. Browning \cite{bro} suggests that one should take
 \begin{equation}
 c_h=\frac{11}{8}f(h)\prod_{p}\left(1-\frac{1}{p}\right)^2\left(1+\frac{2}{p}\right),
 \end{equation}
 for an explicit multiplicative function $f(h)$, and is able to prove that \eqref{D_h} is true on average, namely
 \begin{equation}
 \sum_{h\leq H}\left(D_h(X)-c_h X(\log X)^3\right)=o(HX(\log X)^3)
 \end{equation}
 for $H\geq X^{3/4+\varepsilon}$. 
   
Using spectral theory of automorphic forms, Blomer \cite{blo} improved the range of $H$ substantially to $H\geq X^{1/3+\varepsilon}$ and produced a power saving error term. Furthermore, Blomer's approach seems to be flexible enough to be adapted to the study of more general correlation sums.

\begin{Theorem}[Blomer \cite{blo}\label{blomer's}]
Let $W$ be a smooth function with compact support in $[1,2]$ and Mellin transform $\widehat{W}$. Let $1\leq H\leq X/3$ and $r_d(n)$ be the Ramanujan sum. Let $a(n)$, for $X\leq n\leq 2X$, be any sequence of complex numbers. Then
\begin{equation}
\begin{split}
\sum_{h}W\left(\frac{h}{H}\right)\sum_{X\leq n\leq 2X}a(n)\tau(n+h)\tau(n-h)&=H\widehat{W}(1)\sum_{X\leq n\leq 2X}a(n)\sum_{d}\frac{r_d(2n)}{d^2}(\log n+2\gamma-2\log d)^2
\\&+O\left(X^{\varepsilon}\left(\frac{H^2}{X^{1/2}}+HX^{1/4}+(XH)^{1/2}+\frac{X}{H^{1/2}}\right)\|a\|_2\right).
\end{split}
\end{equation}
\end{Theorem}

Here $\|a\|_2=(\sum|a_n|^2)^{\frac{1}{2}}$ is the $\ell^2$-norm. The first term in the $O$-term above comes from the smooth approximation of the sum on the left. The second one comes from the treatment of the ``minus-case'' after application of Voronoi summation of the divisor function. The last two terms come from spectral methods.

\smallskip
As a corollary, Blomer obtains

\begin{Corollary}[Blomer \cite{blo}]
Let $W$ be a smooth function with compact support in $[1,2]$ and Mellin transform $\widehat{W}$. Let $1\leq H\leq X/3$. Then
\begin{equation}
\begin{split}
\sum_{h}W\left(\frac{h}{H}\right)\sum_{X\leq n\leq 2X}\tau_k(n)\tau(n+h)\tau(n-h)&=\widehat{W}(1)XHQ_{k+1}(\log X)
\\&+O\left(X^{\varepsilon}(H^2+HX^{1-\frac{1}{k+2}}+XH^{1/2}+X^{3/2}H^{-1/2})\right),
\end{split}
\end{equation}
where $\tau_k$ is the $k$-th fold divisor function, and $Q_{k+1}$ is a degree $k+1$ polynomial.
\end{Corollary}

The result is non-trivial for $X^{\frac{1}{3}+\varepsilon}\leq H\leq X^{1-\varepsilon}$, which in the case of $k=2$, substantially improves Browning's result. 
\smallskip

Note that since the divisor function can be viewed as the Fourier coefficient of Eisenstein series, one naturally would ask what will be the case if the divisor functions are replaced by Fourier coefficients of cusp forms. Blomer \cite{blo} remarked that if one uses Jutila's circle method and argues as in \cite{bl-mi}, then one might obtain an analogous result. The purpose of this note is to carry this out in detail. It turns out that new difficulties arise, making it difficult to obtain a range for $H$ as good as the divisor function case. Namely, we are not able to `open' the Fourier coefficient as Blomer does with the divisor function. We obtain the following result.

\begin{Theorem}\label{maintheorem}
Let $1\leq H\leq X/3$. Let $W$ be a smooth function with compact support in $[1,2]$, and $a(n)$, $X\leq n\leq 2X$, be any sequence of complex numbers. Let $\lambda_1(n), \lambda_2(n)$ be Hecke eigenvalues of holomorphic Hecke eigencuspforms of weight $\kappa_1$, $\kappa_2$ for $\mathrm{SL}_2(\mathbb{Z})$ respectively. Then
\begin{equation}\label{finalbdd}
\sum_{h}W\left(\frac{h}{H}\right)\sum_{X\leq n\leq 2X}a(n)\lambda_{1}(n+h)\lambda_{2}(n-h)\ll X^{\varepsilon}\frac{X}{H}\left((XH)^{1/2}+\frac{X}{H^{1/2}}\right)\|a\|_2.
\end{equation}
\end{Theorem}





One should compare our result with the third and fourth terms of the $O$-term in Blomer's theorem, as both of them naturally come from the spectral theory of automorphic forms. The result is nontrivial as long as $H\geq X^{\frac{2}{3}+\varepsilon}$.

As an immediate consequence, we have
\begin{Corollary}
Let $1\leq H\leq X/3$. Let $W$ be a smooth function with compact support in $[1,2]$. Let $\lambda_1(n), \lambda_2(n), \lambda_3(n)$ be Hecke eigenvalues of holomorphic Hecke eigencuspforms of weight $\kappa_1$, $\kappa_2$ and $\kappa_3$ for $\mathrm{SL}_2(\mathbb{Z})$ respectively. Then
\begin{equation}
\sum_{h}W\left(\frac{h}{H}\right)\sum_{X\leq n\leq 2X}\lambda_{1}(n-h)\lambda_{2}(n)\lambda_{3}(n+h)\ll X^{\varepsilon}\min\left(XH,\frac{X^2}{H^{1/2}}\right).
\end{equation}
\end{Corollary}
The result is nontrivial only for $H\geq X^{\frac{2}{3}+\varepsilon}$. One can remove the smooth function $W$ in the $h$-sum, as in \cite{blo}.

\smallskip

If on the other hand, one allows one of the Fourier coefficients to be non-cuspidal, then the advantage to open the divisor function enables us to enlarge the range of $H$ to $H\geq X^{\frac{1}{3}+\varepsilon}$. This feature of decomposable functions was pointed out by Meurman in \cite{Meu}. For instance one has the following result which follows from the same line of proof of Blomer \cite{blo}, although it is not explicitly stated as such in that work.
\begin{Corollary}[Blomer \cite{blo}]
Let $1\leq H\leq X/3$. Let $W$ be a smooth function with compact support in $[1,2]$, and $a(n)$, $X\leq n\leq 2X$, be any sequence of complex numbers. Let $\lambda(n)$ be Hecke eigenvalues of a holomorphic Hecke eigencuspform of weight $\kappa$ for $\mathrm{SL}_2(\mathbb{Z})$. Then
\begin{equation}
\sum_{h}W\left(\frac{h}{H}\right)\sum_{X\leq n\leq 2X}a(n)\tau(n+h)\lambda(n-h)\ll X^{\varepsilon}\left(\frac{H^2}{X^{1/2}}+HX^{1/4}+(XH)^{1/2}+\frac{X}{H^{1/2}}\right)\|a\|_2.
\end{equation}
\end{Corollary}

\smallskip
\bigskip

\textbf{Notation.} $x\asymp X$ means $X\leq x\leq 2X$. We will use the $\varepsilon$-convention: $\varepsilon>0$ is arbitrarily small but not necessarily the same at each occurrence.

\bigskip
\section{Preliminaries}
In this section we collect some lemmas that will be used in our proof.
First let us recall the Voronoi summation formula for holomorphic Hecke eigenvalues.
\begin{Lemma}[{\cite[Theorem A.4]{KMV1}}]\label{voronoi}
Assume $(b,c)=1$. Let $V$ be a smooth compactly supported function. Let $N>0$ and let $\lambda(n)$ denote Hecke eigenvalues of a holomorphic Hecke eigencuspform of weight $\kappa$ for $\mathrm{SL}_2(\mathbb{Z})$. Then
\begin{equation}
\sum_{n}\lambda(n)e\left(\frac{bn}{c}\right)V\left(\frac{n}{N}\right)=\frac{N}{c}\sum_{n}\lambda(n)e\left(-\frac{\bar{b}n}{c}\right)\cdot 2\pi i^\kappa\int_{0}^{\infty}V(x)J_{k-1}\left(\frac{4\pi\sqrt{nNx}}{c}\right)\mathrm{d}x
\end{equation}
\end{Lemma}


\smallskip
We will use the following variant of Jutila's circle method \cite{jut1}, \cite{jut}.
\begin{Lemma}\label{jutila}Let $Q\geq 1$ and $Q^{-2}\leq \delta\leq Q^{-1}$ be two parameters. Let $w$ be a nonnegative function supported in $[Q,2Q]$ and satisfies $\|w\|_\infty\leq1$ and $\sum_{c}w(c)>0$. Let $\mathbb{I}_S$ be the characteristic function of the set $S$. Define\label{jutila}
\begin{equation}
\tilde{I}(\alpha)=\frac{1}{2\delta \Lambda}\sum_{c}w(c)\sumstar_{d(\mathrm{mod} c)}\mathbb{I}_{[\frac{d}{c}-\delta,\frac{d}{c}+\delta]}(\alpha),
\end{equation}
where $\Lambda=\sum_{c}w(c)\phi(c)$. Then $\tilde{I}(\alpha)$ is a good approximation of $\mathbb{I}_{[0,1]}$ in the sense that
\begin{equation}
\int_{0}^{1}|1-\tilde{I}(\alpha)|^2\mathrm{d}\alpha\ll_{\varepsilon} \frac{Q^{2+\varepsilon}}{\delta\Lambda^2}.
\end{equation}
\end{Lemma}

The feature of the circle method in our application, as in \cite{bl-mi},  is that the parameter $Q$ turns out to be just a ``catalyst'', not entering into the final bound. This will become clear at the final stage of our argument.

In order to state Kuznetsov's trace formula and the large sieve inequalities, let us define the following integral transforms for a smooth function $\phi: [0,\infty)\rightarrow \mathbb{C}$ satisfying $\phi(0)=\phi'(0)=0$, $\phi^{(j)}(x)\ll (1+x)^{-3}$ for $0\leq j\leq 3$:
$$\dot{\phi}(k)=4i^k\int_{0}^{\infty}\phi(x)J_{k-1}(x)\frac{\mathrm{d}x}{x},\ \  \tilde{\phi}(t)=2\pi i\int_{0}^{\infty}\phi(x)\frac{J_{2it}(x)-J_{-2it}(x)}{\sinh(\pi t)}\frac{\mathrm{d}x}{x}.$$

We use the notations in \cite{BHM} and \cite{blo}. Let $\mathcal{B}_k$ be an orthonormal basis of the space of holomorphic cusp forms of level $1$ and weight $k$. Let $\mathcal{B}$ be a fixed orthonormal basis of Hecke-Maass eigenforms of level 1. For $f\in \mathcal{B}_k$, we write its Fourier expansion at $\infty$ as
$$f(z)=\sum_{n\geq 1}\rho_f(n)(4\pi n)^{k/2}e(nz).$$

For $f\in \mathcal{B}$ with spectral parameter $t$, we write its Fourier expansion as
$$f(z)=\sum_{n\neq 0}\rho_f(n)W_{0,it}(4\pi|n|y)e(nx),$$
where $W_{0,it}(4\pi|n|y)=(y/\pi)^{1/2}K_{it}(y/2)$ is a Whittaker function. 

Finally for the Eisenstein series $E(z,s)$, we write its Fourier expansion at $s=\frac{1}{2}+it$ as
\begin{equation*}
E\left(z,\frac{1}{2}+it\right)=y^{\frac{1}{2}+it}+\varphi\left(1/2+it\right)y^{\frac{1}{2}-it}+\sum_{n\neq0}\rho(n,t)W_{0,it}(4\pi|n|y)e(nx).
\end{equation*}

Then we have the following Kuznetsov's trace formula in the notation of \cite{BHM}.
\begin{Lemma}\label{kuznetsov} Let $a, b>0$ be integers, then

\begin{equation}
\begin{split}
\sum_{c\geq 1}\frac{1}{c}S(a,b;c)\phi\left(\frac{4\pi\sqrt{ab}}{c}\right)&=\sum_{\substack{k\geq 2\\ k\ \mathrm{even}}}\sum_{f\in\mathcal{B}_k}\dot{\phi}(k)\Gamma(k)\sqrt{ab}\,\overline{\rho_f(a)}\rho_f(b)
\\&+\sum_{f\in\mathcal{B}}\tilde{\phi}(t_f)\frac{\sqrt{ab}}{\cosh (\pi t_f)}\overline{\rho_f(a)}\rho_f(b)+\frac{1}{4\pi}\int_{-\infty}^{\infty}\tilde{\phi}(t)\frac{\sqrt{ab}}{\cosh (\pi t)}\overline{\rho(a,t)}\rho(b,t)\mathrm{d}t.
\end{split}
\end{equation}
\end{Lemma}

We will use the following spectral large sieve inequalities of Deshouillers and Iwaniec \cite{de-iw}.
\begin{Lemma}\label{largesieve} Let $T, M\geq 1$. Let $(a_m)$, $M\leq m\leq 2M$, be a sequence of complex numbers, then all of the three quantities\\
$$\sum_{\substack{2\leq k\leq T\\ k\ \mathrm{even}}}\Gamma(k)\sum_{f\in\mathcal{B}_k}\big|\sum_{m}a_m\sqrt{m}\rho_f(m)\big|^2, \sum_{\substack{f\in\mathcal{B}\\ |t_f|\leq T}}\frac{1}{\cosh (\pi t_f)}\big|\sum_{m}a_m\sqrt{m}\rho_f(\pm m)\big|^2,$$
$$\int_{-T}^{T}\frac{1}{\cosh (\pi t)}\big|\sum_{m}a_m\sqrt{m}\rho(\pm m,t)\big|^2 \mathrm{d}t$$
are bounded by
$$M^{\varepsilon}(T^2+M)\sum_{m}|a_m|^2.$$
\end{Lemma}

We need the following lemma of Blomer-Mili{\'c}evi{\'c} \cite{bl-mi} for a certain type of Bessel transform, which provides the asymptotic behavior of the weight function.
\begin{Lemma}\label{wstar}
Let $W$ be a fixed smooth function with support in $[1,2]$ satisfying $W^{(j)}(x)\ll_j 1$ for all $j$. Let $\nu\in \mathbb{C}$ be a fixed number with $\Re \nu\geq 0$. Define\label{boundofw}
\begin{equation}
W^{\star}(z,w)=\int_{0}^{\infty}W(y)J_{\nu}(4\pi\sqrt{yw+z})\mathrm{d}y.
\end{equation}
Fix $C\geq 1$ and $A,\varepsilon>0$. Then for $z\geq 4|w|>0$ we have
\begin{equation}
W^{\star}(z,w)=W_+(z,w)z^{-1/4}e(2\sqrt{z})+W_-(z,w)z^{-1/4}e(-2\sqrt{z})+O_A(C^{-A})
\end{equation}
for suitable functions $W_\pm$ satisfying
\begin{equation}\label{wplusminus}
z^i w^j\frac{\partial^i}{\partial z^i}\frac{\partial^j}{\partial w^j}W_\pm(z,w)\begin{cases}=0, & \quad   \sqrt{z}/w\leq C^{-\varepsilon}, \\
\ll C^{\varepsilon(i+j)}, & \quad \mathrm{otherwise},\\
\end{cases}
\end{equation}
for any $i,j\in \mathbb{N}_0$. The implied constants depend on $i, j$ and $\nu$.
\end{Lemma}

The lemma holds true for both positive and negative $w$, as long as $z\geq 4|w|$. For a proof, see \cite[Lemma 17]{bl-mi} and the remark after that.

From Lemma \ref{wstar}, one has the following easy consequence. 
\begin{Corollary}[{\cite[Corollary 18]{bl-mi}}]\label{realpart}
Let $\omega$ be a smooth function supported in a rectangle $[c_1,c_2]\times[c_1,c_2]$ for two constants $c_2>c_1>0$, and let $Z\gg 1$, $W>1$ be two parameters such that $c_1Z\geq 4c_2W$, then the double Mellin transform
\begin{equation*}
\widehat{W}_{\pm}(s,t)=\int_{0}^{\infty}\int_{0}^{\infty}W_{\pm}(z,w)\omega\left(\frac{z}{Z},\frac{w}{W}\right)z^sw^t\frac{\mathrm{d}z}{z}\frac{\mathrm{d}w}{w}
\end{equation*}
is holomorphic on $\mathbb{C}^2$, and satisfies
\begin{equation*}
\widehat{W}_{\pm}(s,t)\ll_{A,B,\varepsilon,\Re s, \Re t} C^{\varepsilon}(1+|s|)^{-A}(1+|t|)^{-B}
\end{equation*}
on vertical lines. 
\end{Corollary}
We will use this corollary when we try to separate variables, after the application of Kuznetsov's trace formula. Then, the following lemma of \cite{blo} will help us to truncate the lengths of summations. See also \cite[Lemma 3]{jut}.
\begin{Lemma}[{\cite[Lemma 5]{blo}}] \label{lemmaofphi}
Let $\mathcal{Z}\gg1$,$\tau\in\mathbb{R}$, $\alpha\in[-4/5,4/5]$ and $w$ be a smooth compactly supported function. For
\begin{equation}
\phi(z)=e^{\pm iz\alpha}w\left(\frac{z}{\mathcal{Z}}\right)\left(\frac{z}{\mathcal{Z}}\right)^{i\tau},
\end{equation}
we have
\begin{equation}
\dot{\phi}(k)\ll_A\frac{1+|\tau|}{\mathcal{Z}}\left(1+\frac{k}{\mathcal{Z}}\right)^{-A}, \tilde{\phi}(t)\ll_A\left(1+\frac{|t|+\mathcal{Z}}{1+|\tau|}\right)^{-A}
\end{equation}
for $t\in\mathbb{R}$, $k\in\mathbb{N}$ and any $A\geq 0$.
\end{Lemma}
\begin{Remark} From Lemma \ref{lemmaofphi} we see that $\dot{\phi}(k)$ is negligibly small as long as $|k|>\mathcal{Z}$, so that later one can truncate summation over $k$ at $\mathcal{Z}$, up to a negligible error. Moreover, in application, $\mathcal{Z}$ will usually be relatively larger than $\tau$, so that the contribution from $\tilde{\phi}(t)$ is always negligible. In particular, after the application of Kuznetsov's trace formula later, we only need to treat the contribution from the holomorphic spectrum, since the Maass forms and Eisenstein series parts will be negligible due to the rapid decay of the weight function $\tilde{\phi}(t)$. 
\end{Remark}

\bigskip
\section{Proof of the theorem}
Let $\lambda_1(n), \lambda_2(n)$ be the Hecke eigenvalues of holomorphic Hecke eigencuspforms of weight $\kappa_1$, $\kappa_2$ for $\mathrm{SL}_2(\mathbb{Z})$. We change the order of summation of $n$ and $h$, fix $n\asymp X$, and first deal with the sum over $h$
\begin{equation}\label{originalsum}
E(n):=\sum_{h\asymp H}\lambda_{1}\left(n+h\right)\lambda_{2}\left(n-h\right)W\left(\frac{h}{H}\right).
\end{equation}

Let $n+h=m_1$ and $n-h=m_2$, then $m_1+m_2=2n$, and we can rewrite the summation above as

\begin{equation}\label{shiftedsum}
\sum_{m_1+m_2=2n}\lambda_{1}(m_1)\lambda_{2}(m_2)W\left(\frac{m_1-n}{H}\right)V\left(\frac{n-m_2}{H'}\right).
\end{equation}
Here $V$ is a redundant smooth function used to keep track of the support of $m_2$, and $H'$ is a parameter to be determined later, satisfying $H\leq H'\leq X/3$. Later we will see that in our case, $H'=X/3$ will give the best result. Here one can take, say, $V=W$.
\smallskip

\begin{Remark}  Let us make a comment on the difference between the shifted convolution sum in \eqref{shiftedsum} with the one treated by Blomer and Mili{\'c}evi{\'c}, \cite[(3.7)]{bl-mi}. In our case, we have localized both the variables $m_1$ and $m_2$ to vary in intervals \textbf{around} $n$. In Blomer and Mili{\'c}evi{\'c}'s case, one of the variables varies in an interval of length $2H$ around $n$, while the other variable varies, say, in $[H/2,2H]$. One should also note that for $f$ a holomorphic cusp form, the sum $\sum_{m_1+m_2=X}\lambda_{f}(m_1)\lambda_{f}(m_2)$ is just the $X$-th Fourier coefficient of the cusp form $f^2$ and thus one can apply Deligne's estimate. 
\end{Remark}

Now we follow the approach of Jutila \cite{jut} and Blomer and Mili{\'c}evi{\'c} \cite[sections 7 $\&$ 8]{bl-mi}. Let $C=X^{1000}$ be a large parameter. Apply Lemma \ref{jutila} with $Q=C$ and $\delta=C^{-1}$. Let $w_0$ be a fixed smooth function supported in $[1,2]$. Let $w(c)=w_0(c/C)$. In particular we have $\Lambda=\sum_{c}w_0(c/C)\phi(c)\gg C^{2-\varepsilon}$. Detecting the condition $m_1+m_2=2n$ by $\int_{0}^{1}e(\alpha(m_1+m_2-2n))\mathrm{d} \alpha$ and applying Jutila's circle method, we have
\begin{equation*}
E(n)=\widetilde{E}(n)+\mathrm{Error},
\end{equation*}
where
\begin{equation}
\begin{split}
\widetilde{E}(n)=\frac{1}{2\delta}\int_{-\delta}^{\delta}&\frac{1}{\Lambda}\sum_{c}w_0\left(\frac{c}{C}\right)\sumstar_{d(c)}e\left(\frac{-2nd}{c}\right)\sum_{m_1}\lambda_{1}(m_1)e\left(\frac{dm_1}{c}\right)W\left(\frac{m_1-n}{H}\right)e(\eta(m_1-n))
\\&\cdot \sum_{m_2}\lambda_{2}(m_2)e\left(\frac{dm_2}{c}\right)V\left(\frac{n-m_2}{H'}\right)e(-\eta(n-m_2))\mathrm{d}\eta,
\end{split}
\end{equation}
and
\begin{equation*}
\begin{split}
\mathrm{Error}&=\int_{0}^{1}\sum_{m_1}\sum_{m_2}\lambda_{1}(m_1)\lambda_{2}(m_2)W\left(\frac{m_1-n}{H}\right)V\left(\frac{n-m_2}{H'}\right)e(\alpha(m_1+m_2-2n))(1-\tilde{I}(\alpha))\mathrm{d} \alpha\\
&\ll\left(\sum_{m_1}\left|\lambda_{1}(m_1)W\big(\frac{m_1-n}{H}\big)\right|\right)\left(\sum_{m_2}\left|\lambda_{2}(m_2)V\big(\frac{n-m_2}{H'}\big)\right|\right)\cdot \frac{C^{1+\varepsilon}}{\delta^{1/2}\Lambda}\\
&\ll\frac{X^{2}C^{1+\varepsilon}}{\delta^{1/2}\Lambda}\ll C^{-2/5}
\end{split}
\end{equation*}
by the Cauchy-Schwarz inequality, Lemma \ref{jutila} and the bound $\sum_{n\leq x}|\lambda_f(n)|^2\ll_f x$.
\smallskip

Denote $W_\eta(x)=W(x)e(\eta x)$ and $V_{-\eta}(x)=V(x)e(-\eta x)$. Then
\begin{equation}
\begin{split}
\widetilde{E}(n)=\frac{1}{2\delta}\int_{-\delta}^{\delta}&\frac{1}{\Lambda}\sum_{c}w_0\left(\frac{c}{C}\right)\sumstar_{d(c)}e\left(\frac{-2nd}{c}\right)\sum_{m_1}\lambda_{1}(m_1)e\left(\frac{dm_1}{c}\right)W_{\eta H}\left(\frac{m_1-n}{H}\right)\\
&\cdot \sum_{m_2}\lambda_{2}(m_2)e\left(\frac{dm_2}{c}\right)V_{-\eta H'}\left(\frac{n-m_2}{H'}\right)\mathrm{d}\eta.
\end{split}
\end{equation}

Note that since $|\eta|\leq C^{-1}=X^{-1000}$ is very small, the functions $W_{\eta H}$ and $V_{-\eta H'}$  have nice properties inherited from $W$. In particular, one has $W^{(j)}_{\eta H},\ V^{(j)}_{-\eta H'}\ll 1$ for any $j\geq 0$, supported in $[1,2]$, uniformly for $|\eta|\leq C^{-1}$. 

Applying Voronoi summation formula to the $m_1$-sum, we have
\begin{equation}
\sum_{m_1}\lambda_{1}(m_1)e\left(\frac{dm_1}{c}\right)W_{\eta H}\left(\frac{m_1-n}{H}\right)=\frac{H}{c}\sum_{m_1}\lambda_{1}(m_1)e\left(-\frac{\bar{d}m_1}{c}\right)W^{\star}_{\eta H}\left(\frac{m_1n}{c^2},\frac{m_1H}{c^2}\right),
\end{equation}
where
\begin{equation}
W^{\star}_{\eta H}(z,w)=2\pi i^{\kappa_1}\int_{0}^{\infty}W_{\eta H}(y)J_{\kappa_1-1}(4\pi\sqrt{yw+z})\mathrm{d}y.
\end{equation}

Similarly,
\begin{equation}
\sum_{m_2}\lambda_{2}(m_2)e\left(\frac{dm_2}{c}\right)V_{-\eta H'}\left(\frac{n-m_2}{H'}\right)=\frac{H'}{c}\sum_{m_2}\lambda_{2}(m_2)e\left(-\frac{\bar{d}m_2}{c}\right)V^{\star}_{-\eta H'}\left(\frac{m_2n}{c^2},\frac{m_2H'}{c^2}\right),
\end{equation}
with
\begin{equation}
V^{\star}_{-\eta H'}(z,w)=2\pi i^{\kappa_2}\int_{0}^{\infty}V_{-\eta H'}(y)J_{\kappa_2-1}(4\pi\sqrt{-yw+z})\mathrm{d}y.
\end{equation}

Substituting these back into $\widetilde{E}(n)$, we get

\begin{equation}
\widetilde{E}(n)=\frac{1}{2\delta}\int_{-\delta}^{\delta} \widetilde{E}_\eta(n)\mathrm{d}\eta,
\end{equation}
with
\begin{equation}
\begin{split}
\widetilde{E}_\eta(n)=&\frac{HH'}{\Lambda}\sum_{c}\frac{w_0(c/C)}{c^2}\sum_{m_1}\sum_{m_2}\lambda_{1}(m_1)\lambda_{2}(m_2)S(m_1+m_2,2n;c)\\
&\cdot W^{\star}_{\eta H}\left(\frac{m_1n}{c^2},\frac{m_1H}{c^2}\right)V^{\star}_{-\eta H'}\left(\frac{m_2n}{c^2},\frac{m_2H'}{c^2}\right).
\end{split}
\end{equation}

If one can establish a bound for $\sum_{n\asymp X}a(n)\widetilde{E}_\eta(n)$, uniformly in $\eta$, then the same bound will hold true for $\sum_{n\asymp X}a(n)E(n)$, and we are done.

By Lemma \ref{boundofw}, we can write
\begin{equation*}
W^{\star}_{\eta H}\left(\frac{m_1n}{c^2},\frac{m_1H}{c^2}\right)=\sum_{\pm}W_{\pm}\left(\frac{m_1n}{c^2},\frac{m_1H}{c^2}\right)\left(\frac{m_1n}{c^2}\right)^{-\frac{1}{4}}e\left(\pm2\frac{\sqrt{m_1n}}{c}\right)+O\left(C^{-A}\right),
\end{equation*}
for some $W_{\pm}$ satisfying \eqref{wplusminus}.

Similarly we have
\begin{equation*}V^{\star}_{-\eta H'}\left(\frac{m_2n}{c^2},\frac{m_2H'}{c^2}\right)=\sum^{\pm}V^{\pm}\left(\frac{m_2n}{c^2},\frac{m_2H'}{c^2}\right)\left(\frac{m_2n}{c^2}\right)^{-\frac{1}{4}}e\left(\pm2\frac{\sqrt{m_2n}}{c}\right)+O\left(C^{-A}\right),
\end{equation*}
for some $V^{\pm}$ satisfying similar conditions.
Note that from now on we will not display dependence on the parameter $\eta$ anymore, since all of the estimations we obtain will be uniform in $\eta$, by our previous remark on the nice properties of $W_{\eta H}$ and $V_{-\eta H'}$.

By the bound \eqref{wplusminus}, we can restrict the lengths of the new sums over $m_1$, $m_2$ to
\begin{equation*}
m_1\leq T_1:=\frac{C^{2+\varepsilon}X}{H^2},\ m_2\leq T_2:=\frac{C^{2+\varepsilon}X}{H'^2},
\end{equation*}
up to a negligible error.

Now we further restrict the lengths of summations to dyadic segments. That is, we assume

$$m_1\asymp \mathcal{M}_1,\ m_2\asymp \mathcal{M}_2,$$
where $\mathcal{M}_1\ll T_1,\ \mathcal{M}_2\ll T_2$.

Denote $b=m_1+m_2$. Then
$$b\asymp \mathcal{M}_1+\mathcal{M}_2.$$

Then $\widetilde{E}(n)$ will be a sum of at most $O((\log C)^3)$ terms, of the following form.
\begin{equation}
\begin{split}
\widetilde{E}(n,\mathcal{M}_1,\mathcal{M}_2)&=\frac{HH'n^{-\frac{1}{2}}}{\Lambda}\sum_{b\asymp \mathcal{M}_1+\mathcal{M}_2}\sum_{\substack{m_1+m_2=b\\ m_1\asymp\mathcal{M}_1,\ m_2\asymp\mathcal{M}_2}}\lambda_{1}(m_1)m_1^{-\frac{1}{4}}\lambda_{2}(m_2)m_2^{-\frac{1}{4}}\sum_{c}\frac{S(b,2n;c)}{c}w_0\left(\frac{c}{C}\right)
\\&
\cdot\sum_{\pm}\sum^{\pm}W_{\pm}\left(\frac{m_1n}{c^2},\frac{m_1H}{c^2}\right)V^{\pm}\left(\frac{m_2n}{c^2},\frac{m_2H'}{c^2}\right)e\left(\pm2\frac{\sqrt{m_1n}}{c}\right)e\left(\pm2\frac{\sqrt{m_2n}}{c}\right)+O(C^{-10}).
\end{split}
\end{equation}

To prepare for the application of Kuznetsov's trace formula, we combine the weight functions above and write
\begin{equation}
\widetilde{E}(n,\mathcal{M}_1,\mathcal{M}_2)=\frac{HH'n^{-\frac{1}{2}}}{\Lambda}\sum_{b\asymp \mathcal{M}_1+\mathcal{M}_2}\sum_{\substack{m_1+m_2=b\\ m_1\asymp\mathcal{M}_1,\ m_2\asymp\mathcal{M}_2}}\lambda_{1}(m_1)m_1^{-\frac{1}{4}}\lambda_{2}(m_2)m_2^{-\frac{1}{4}}\sum_{c}\frac{S(b,2n;c)}{c}\Phi\left(\frac{4\pi\sqrt{2nb}}{c}\right),
\end{equation}
where
\begin{equation}
\Phi(z)=\sum_{\pm}\sum^{\pm}w_0\left(\frac{4\pi\sqrt{2nb}}{zC}\right)W_{\pm}\left(\frac{m_1z^2}{32\pi^2b},\frac{m_1Hz^2}{32\pi^2nb}\right)
V^{\pm}\left(\frac{m_2z^2}{32\pi^2b},\frac{m_2H'z^2}{32\pi^2nb}\right)e^{\pm i\sqrt{\frac{m_1}{2b}}z}e^{\pm i\sqrt{\frac{m_2}{2b}}z}.
\end{equation}

Note that the support of $w_0$ implies that we can restrict $z$ to
$$z\asymp \mathcal{Z}:=\frac{\sqrt{X(\mathcal{M}_1+\mathcal{M}_2)}}{C}.$$
We can attach a redundant smooth weight function $w_2(\frac{z}{\mathcal{Z}})$ of compact support $[1/100\mathcal{Z},100\mathcal{Z}]$ that is constantly $1$ on $[1/20\mathcal{Z},20\mathcal{Z}]$, to $\Phi(z)$.

Now we separate the variables. We do this by Mellin inversion to the functions $w_0$, $W_{\pm}$ and $V^{\pm}$ (using Corollary \ref{realpart}). This can be done with almost no loss since these functions are non-oscillatory, similar to \cite{bl-mi} and \cite{bfkmm}. Thus we have
\begin{equation}
\begin{split}
\Phi(z)=&\sum_{\pm}\sum^{\pm}w_0\left(\frac{4\pi\sqrt{2nb}}{zC}\right)W_{\pm}\left(\frac{m_1z^2}{32\pi^2b},\frac{m_1Hz^2}{32\pi^2nb}\right)V^{\pm}\left(\frac{m_2z^2}{32\pi^2b},\frac{m_2H'z^2}{32\pi^2nb}\right)\\
&\cdot e^{\pm iz(\sqrt{\frac{m_1}{2b}}\pm\sqrt{\frac{m_2}{2b}})}w_2\left(\frac{z}{\mathcal{Z}}\right)\\
=&\sum_{\pm}\sum^{\pm}\int_{\mathcal{C}}\widehat{w}_0(s_1)\widehat{W}_{\pm}(s_2,s_3)\widehat{V}^{\pm}(s_4, s_5)
\bigg(\frac{4\pi\sqrt{2nb}}{zC}\bigg)^{-s_1}\bigg(\frac{m_1z^2}{32\pi^2b}\bigg)^{-s_2}\\
&\cdot\bigg(\frac{m_1Hz^2}{32\pi^2nb}\bigg)^{-s_3}\bigg(\frac{m_2z^2}{32\pi^2b}\bigg)^{-s_4}\bigg(\frac{m_2H'z^2}{32\pi^2nb}\bigg)^{-s_5}\mathrm{d} s\cdot e^{\pm iz(\sqrt{\frac{m_1}{2b}}\pm\sqrt{\frac{m_2}{2b}})}w_2\left(\frac{z}{\mathcal{Z}}\right)\\
=&\sum_{\pm}\sum^{\pm}\int_{\mathcal{C}}\widehat{w}_0(s_1)\widehat{W}_{\pm}(s_2,s_3)\widehat{V}^{\pm}(s_4, s_5)(4\pi\sqrt{2})^{-s_1}(32\pi^2)^{s_2+s_3+s_4+s_5}C^{s_1}\mathcal{Z}^{s_1-2s_2-2s_3-2s_4-2s_5}\\
&\cdot (\mathcal{M}_1+\mathcal{M}_2)^{-\frac{s_1}{2}+s_2+s_3+s_4+s_5}\mathcal{M}_1^{-s_2-s_3}\mathcal{M}_2^{-s_4-s_5}H^{-s_3}H'^{-s_5}n^{-\frac{s_1}{2}+s_3+s_5}\left(\frac{b}{\mathcal{M}_1+\mathcal{M}_2}\right)^{-\frac{s_1}{2}+s_2+s_3+s_4+s_5}\\
&\cdot\bigg(\frac{m_1}{\mathcal{M}_1}\bigg)^{-s_2-s_3}\bigg(\frac{m_2}{\mathcal{M}_2}\bigg)^{-s_4-s_5}
\bigg(\frac{z}{\mathcal{Z}}\bigg)^{s_1-2s_2-2s_3-2s_4-2s_5}e^{\pm iz(\sqrt{\frac{m_1}{2b}}\pm\sqrt{\frac{m_2}{2b}})}w_2\left(\frac{z}{\mathcal{Z}}\right)\mathrm{d}s,
\end{split}
\end{equation}
where $\mathcal{C}$ is the fivefold contour taken over the vertical lines $\Re s_1=\Re s_2=\Re s_3=\Re s_4=\Re s_5=0$.
Here we denote $\mathrm{d}s=\frac{1}{(2\pi i)^5}\prod_{j=1}^{5}\mathrm{d} s_j$. Note that due to the rapid decay of $\widehat{w}_1(s_1)\widehat{W}_{\pm}(s_2,s_3)\widehat{V}^{\pm}(s_4, s_5)$ along vertical lines, we can truncate the integrals above at $|\Im s_i|\leq C^\varepsilon, 1\leq i\leq 5$, at the cost of a negligible error. We denote the truncated contour by $\widetilde{\mathcal{C}}$.

We arrive at
\begin{equation}
\begin{split}
\widetilde{E}(n,\mathcal{M}_1,\mathcal{M}_2)&=\frac{HH'}{\Lambda}(\mathcal{M}_1\mathcal{M}_2)^{-\frac{1}{4}}\sum_{\pm}\sum^{\pm}\int_{\widetilde{\mathcal{C}}}\widehat{w}_0(s_1)\widehat{W}_{\pm}(s_2,s_3)\widehat{V}^{\pm}(s_4, s_5)(4\pi\sqrt{2})^{-s_1}(32\pi^2)^{s_2+s_3+s_4+s_5}
\\&\cdot C^{s_1}\mathcal{Z}^{s_1-2s_2-2s_3-2s_4-2s_5}(\mathcal{M}_1+\mathcal{M}_2)^{-\frac{s_1}{2}+s_2+s_3+s_4+s_5}\mathcal{M}_1^{-s_2-s_3}\mathcal{M}_2^{-s_4-s_5}H^{-s_3}H'^{-s_5}\\
&\cdot n^{-\frac{1}{2}-\frac{s_1}{2}+s_3+s_5}\sum_{b\asymp \mathcal{M}_1+\mathcal{M}_2}\left(\frac{b}{\mathcal{M}_1+\mathcal{M}_2}\right)^{-\frac{s_1}{2}+s_2+s_3+s_4+s_5}\sum_{\substack{m_1+m_2=b\\ m_1\asymp\mathcal{M}_1,\ m_2\asymp\mathcal{M}_2}}\lambda_{1}(m_1)\lambda_{2}(m_2)
\\&\cdot\left(\frac{m_1}{\mathcal{M}_1}\right)^{-\frac{1}{4}-s_2-s_3}\left(\frac{m_2}{\mathcal{M}_2}\right)^{-\frac{1}{4}-s_4-s_5}\sum_{c}\frac{S(b,2n;c)}{c}\Theta\left(\frac{4\pi\sqrt{2nb}}{c}\right)\mathrm{d}s+O(C^{-10}),
\end{split}
\end{equation}
where
\begin{equation}\label{function}
\Theta(z)=e^{\pm iz(\sqrt{\frac{m_1}{2b}}\pm\sqrt{\frac{m_2}{2b}})}w_2\left(\frac{z}{\mathcal{Z}}\right)\left(\frac{z}{\mathcal{Z}}\right)^{s_1-2s_2-2s_3-2s_4-2s_5}.
\end{equation}

\smallskip

Now we apply Kuznetsov's trace formula to the $c$-sum. By Lemma \ref{lemmaofphi}, the spectral sums can be truncated at $C^{\varepsilon}\mathcal{Z}$. Hence we obtain
\begin{equation}
\sum_{c}\frac{S(b,2n;c)}{c}\Theta\left(\frac{4\pi\sqrt{2nb}}{c}\right)=\mathcal{H}(n)+\mathcal{M}(n)+\mathcal{E}(n)+\text{(negligible error)},
\end{equation}
where
\begin{equation}\label{ksum}
\begin{split}
& \mathcal{H}(n)=\sum_{\substack{2\leq k\leq C^{\varepsilon}\mathcal{Z}\\ k\ \mathrm{even}}}\sum_{f\in\mathcal{B}_k}\Gamma(k)\cdot4i^k\int_{0}^{\infty}\Theta(z)J_{k-1}(z)\frac{\mathrm{d}z}{z}\cdot \sqrt{2n}\rho_f(2n)\sqrt{b}\,\overline{\rho_f(b)}, \\
& \mathcal{M}(n)=\sum_{\substack{f\in\mathcal{B}\\|t_f|\leq C^{\varepsilon}\mathcal{Z}}}2\pi i\int_{0}^{\infty}\Theta(z)\frac{J_{2it_f}(z)-J_{-2it_f}(z)}{\sinh(\pi t_f)}\frac{\mathrm{d}z}{z}\cdot \frac{\sqrt{2nb}}{\cosh (\pi t_f)}\rho_f(2n)\overline{\rho_f(b)}, \\
& \mathcal{E}(n)=\frac{1}{4\pi}\int_{|t|\leq C^{\varepsilon}\mathcal{Z}}2\pi i\int_{0}^{\infty}\Theta(z)\frac{J_{2it}(z)-J_{-2it}(z)}{\sinh(\pi t)}\frac{\mathrm{d}z}{z}\cdot \frac{\sqrt{2nb}}{\cosh (\pi t)}\rho(2n,t)\overline{\rho(b,t)}\mathrm{d}t
\end{split}
\end{equation}
are contributions from the holomorphic modular forms, Maass forms, and Eisenstein series respectively.
\smallskip

\smallskip

Now we are ready to sum over $a(n)$, $X\leq n\leq 2X$. Summing over $n$, we get

\begin{equation}
\begin{split}
&\sum_{n\asymp X}a(n)\widetilde{E}(n,\mathcal{M}_1,\mathcal{M}_2)\\
&=\frac{HH'}{\Lambda}(\mathcal{M}_1\mathcal{M}_2)^{-\frac{1}{4}}\sum_{\pm}\sum^{\pm}\int_{\widetilde{\mathcal{C}}}\widehat{w}_0(s_1)\widehat{W}_{\pm}(s_2,s_3)\widehat{V}^{\pm}(s_4, s_5)(4\pi\sqrt{2})^{-s_1}(32\pi^2)^{s_2+s_3+s_4+s_5}\\
&\cdot C^{s_1}\mathcal{Z}^{s_1-2s_2-2s_3-2s_4-2s_5}(\mathcal{M}_1+\mathcal{M}_2)^{-\frac{s_1}{2}+s_2+s_3+s_4+s_5}\mathcal{M}_1^{-s_2-s_3}\mathcal{M}_2^{-s_4-s_5}H^{-s_3}H'^{-s_5}\\
&\cdot\sum_{n\asymp X}a(n)n^{-\frac{1}{2}-\frac{s_1}{2}+s_3+s_5}\sum_{b\asymp \mathcal{M}_1+\mathcal{M}_2}\left(\frac{b}{\mathcal{M}_1+\mathcal{M}_2}\right)^{-\frac{s_1}{2}+s_2+s_3+s_4+s_5}\sum_{\substack{m_1+m_2=b\\ m_1\asymp\mathcal{M}_1,\ m_2\asymp\mathcal{M}_2}}\lambda_{1}(m_1)\lambda_{2}(m_2)\\
&\cdot\bigg(\frac{m_1}{\mathcal{M}_1}\bigg)^{-\frac{1}{4}-s_2-s_3}\bigg(\frac{m_2}{\mathcal{M}_2}\bigg)^{-\frac{1}{4}-s_4-s_5}\bigg(\mathcal{H}(n)+\mathcal{M}(n)+\mathcal{E}(n)\bigg)\mathrm{d}s+O(C^{-10})\\
&:=\mathcal{HH}+\mathcal{MM}+\mathcal{EE}+O(C^{-10}),
\end{split}
\end{equation}
with $\mathcal{HH}$, $\mathcal{MM}$, $\mathcal{EE}$ being contributions from holomorphic forms, Maass forms and Eisenstein series respectively.
\smallskip

Note that for the function $\Theta(z)$ in \eqref{function}, the imaginary part $\Im (s_1-2s_2-2s_3-2s_4-2s_5)\ll C^{\varepsilon}$, which is relatively small compared to $\mathcal{Z}$. By Lemma \ref{lemmaofphi} and the remark after it, it suffices to deal with the contribution from the holomorphic spectrum, since the contributions from the other two parts will be similar or even smaller.
Also note that since we are considering the full modular group case, we do not have exceptional eigenvalues contribution in the Maass spectrum.

Now we focus on the holomorphic contribution, which is
\begin{equation}
\begin{split}
\mathcal{HH}&=\frac{C^{\varepsilon}HH'}{\Lambda}(\mathcal{M}_1\mathcal{M}_2)^{-\frac{1}{4}}\sum_{\pm}\sum^{\pm}\int_{\widetilde{\mathcal{C}}}\widehat{w}_0(s_1)\widehat{W}_{\pm}(s_2,s_3)\widehat{V}^{\pm}(s_4, s_5)(4\pi\sqrt{2})^{-s_1}(32\pi^2)^{s_2+s_3+s_4+s_5}
\\&\cdot C^{s_1}\mathcal{Z}^{s_1-2s_2-2s_3-2s_4-2s_5}(\mathcal{M}_1+\mathcal{M}_2)^{-\frac{s_1}{2}+s_2+s_3+s_4+s_5}\mathcal{M}_1^{-s_2-s_3}\mathcal{M}_2^{-s_4-s_5}H^{-s_3}H'^{-s_5}\\
&\cdot\int_{0}^{\infty}\sum_{\substack{2\leq k\leq C^{\varepsilon}\mathcal{Z}\\ k\ \mathrm{even}}}4i^k\Gamma(k)\sum_{f\in\mathcal{B}_k}\sum_{n\asymp X}a(n)n^{-\frac{1}{2}-\frac{s_1}{2}+s_3+s_5}\sqrt{2n}\rho_f(2n)\cdot \sum_{b\asymp \mathcal{M}_1+\mathcal{M}_2}\sqrt{b}\,\overline{\rho_f(b)}\\
&\cdot\left(\frac{b}{\mathcal{M}_1+\mathcal{M}_2}\right)^{-\frac{s_1}{2}+s_2+s_3+s_4+s_5}
\sum_{\substack{m_1+m_2=b\\ m_1\asymp\mathcal{M}_1,\ m_2\asymp\mathcal{M}_2}}\lambda_{1}(m_1)\lambda_{2}(m_2)\left(\frac{m_1}{\mathcal{M}_1}\right)^{-\frac{1}{4}-s_2-s_3}
\left(\frac{m_2}{\mathcal{M}_2}\right)^{-\frac{1}{4}-s_4-s_5}\\
&\cdot e^{\pm iz(\sqrt{m_1}\pm\sqrt{m_2})}w_2\left(\frac{\sqrt{2b}z}{\mathcal{Z}}\right)\left(\frac{\sqrt{2b}z}{\mathcal{Z}}\right)^{s_1-2s_2-2s_3-2s_4-2s_5}J_{k-1}(\sqrt{2b}z)\frac{\mathrm{d}z}{z}\mathrm{d}s,
\end{split}
\end{equation}
after making the change of variable $z\mapsto \sqrt{2b}z$. The support of the smooth function $w_2$ implies that $z\asymp \frac{\sqrt{X}}{C}$ in the inner integral.

Bounding the $z$-integral and $s_i$-integrals trivially, we have
\begin{equation}
\begin{split}
\mathcal{HH}&\ll\frac{C^{\varepsilon}HH'}{\Lambda}(\mathcal{M}_1\mathcal{M}_2)^{-\frac{1}{4}}\cdot
\sup_{\substack{|u_1|, |u_2|,|u_3|,|u_4|\leq C^{\varepsilon}\\z\asymp \frac{\sqrt{X}}{C}}}\bigg|\sum_{\substack{2\leq k\leq C^{\varepsilon}\mathcal{Z}\\ k\ \mathrm{even}}}\Gamma(k)\sum_{f\in\mathcal{B}_k}\sum_{n\asymp X}a(n)n^{-\frac{1}{2}+iu_4}\sqrt{2n}\rho_f(2n)\\
&\cdot \sum_{b\asymp \mathcal{M}_1+\mathcal{M}_2}\sqrt{b}\,\overline{\rho_f(b)}J_{k-1}(\sqrt{2b}z)\cdot w_2\left(\frac{\sqrt{2b}z}{\mathcal{Z}}\right)\left(\frac{\sqrt{2b}z}{\mathcal{Z}}\right)^{-2iu_3}\left(\frac{b}{\mathcal{M}_1+\mathcal{M}_2}\right)^{iu_3}\\
&\cdot\sum_{\substack{m_1+m_2=b\\ m_1\asymp\mathcal{M}_1,\ m_2\asymp\mathcal{M}_2}}\lambda_{1}(m_1)\lambda_{2}(m_2)\left(\frac{m_1}{\mathcal{M}_1}\right)^{-\frac{1}{4}+iu_1}\left(\frac{m_2}{\mathcal{M}_2}\right)^{-\frac{1}{4}+iu_2}e^{\pm iz(\sqrt{m_1}\pm\sqrt{m_2})}\bigg|.
\end{split}
\end{equation}

The Cauchy-Schwarz inequality further implies that
\begin{equation}\label{twosup}
\begin{split}
\mathcal{HH}&\ll\frac{C^{\varepsilon}HH'}{\Lambda}(\mathcal{M}_1\mathcal{M}_2)^{-\frac{1}{4}}\cdot
\bigg(\sup_{|u_4|\leq C^{\varepsilon}}\sum_{\substack{2\leq k\leq C^{\varepsilon}\mathcal{Z}\\ k\ \mathrm{even}}}\Gamma(k)\sum_{f\in\mathcal{B}_k}\bigg |\sum_{n\asymp X}a(n)n^{-\frac{1}{2}+iu_4}\sqrt{2n}\rho_f(2n)\bigg|^2\bigg)^{\frac{1}{2}}
\\&\cdot\bigg(\sup_{\substack{|u_1|, |u_2|,|u_3|\leq C^{\varepsilon}\\z\asymp \frac{\sqrt{X}}{C}}}\sum_{\substack{2\leq k\leq C^{\varepsilon}\mathcal{Z}\\ k\ \mathrm{even}}}\Gamma(k)\sum_{f\in\mathcal{B}_k}\bigg|\sum_{b\asymp \mathcal{M}_1+\mathcal{M}_2}\sqrt{b}\,\overline{\rho_f(b)}J_{k-1}(\sqrt{2b}z)\gamma^{\star}(b,z)\bigg|^2\bigg)^{\frac{1}{2}},
\end{split}
\end{equation}
where
\begin{equation}\label{gammastar}
\begin{split}
\gamma^{\star}(b,z)&:=w_2\left(\frac{\sqrt{2b}z}{\mathcal{Z}}\right)\left(\frac{\sqrt{2b}z}{\mathcal{Z}}\right)^{-2iu_3}\left(\frac{b}{\mathcal{M}_1+\mathcal{M}_2}\right)^{iu_3}\\
&\cdot\sum_{\substack{m_1+m_2=b\\ m_1\asymp\mathcal{M}_1,\ m_2\asymp\mathcal{M}_2}}\lambda_{1}(m_1)\lambda_{2}(m_2)\left(\frac{m_1}{\mathcal{M}_1}\right)^{-\frac{1}{4}+iu_1}e^{\pm iz\sqrt{m_1}}\left(\frac{m_2}{\mathcal{M}_2}\right)^{-\frac{1}{4}+iu_2}e^{\pm iz\sqrt{m_2}}.
\end{split}
\end{equation}

The large sieve inequalities yield that
\begin{equation*}
\bigg(\sup_{|u_4|\leq C^{\varepsilon}}\sum_{\substack{2\leq k\leq C^{\varepsilon}\mathcal{Z}\\ k\ \mathrm{even}}}\Gamma(k)\sum_{f\in\mathcal{B}_k}\bigg |\sum_{n\asymp X}a(n)n^{-\frac{1}{2}+iu_4}\sqrt{2n}\rho_f(2n)\bigg|^2\bigg)^{\frac{1}{2}}
\ll \mathcal{C}^{\varepsilon}\left(\mathcal{Z}^2+X\right)^{\frac{1}{2}}X^{-\frac{1}{2}}\|a\|_2.
\end{equation*}

Now it remains to deal with the second line of \eqref{twosup}. We denote
\begin{equation*}
(\star\star):=\sup_{\substack{|u_1|, |u_2|,|u_3|\leq C^{\varepsilon}\\z\asymp \frac{\sqrt{X}}{C}}}\sum_{\substack{2\leq k\leq C^{\varepsilon}\mathcal{Z}\\ k\ \mathrm{even}}}\Gamma(k)\sum_{f\in\mathcal{B}_k}\bigg|\sum_{b\asymp \mathcal{M}_1+\mathcal{M}_2}\sqrt{b}\,\overline{\rho_f(b)}J_{k-1}(\sqrt{2b}z)\gamma^{\star}(b,z)\bigg|^2.
\end{equation*}
We want to separate the $k$-variable from the argument of the $J$-Bessel function, so that one can apply the large sieve inequalities. By the integral representation $J_{k-1}(x)=\frac{1}{\pi}\int_0^{\pi}\cos((k-1)\xi-x\sin\xi)\mathrm{d}\xi$, \cite[(8.411.1)]{GR07}, we have
\begin{equation*}
J_{k-1}(\sqrt{2b}z)=\frac{1}{2\pi}\int_0^{\pi}\left(e^{i(k-1)\xi}e^{-i\sqrt{2b}z\sin\xi}+e^{-i(k-1)\xi}e^{i\sqrt{2b}z\sin\xi}\right)\mathrm{d}\xi.
\end{equation*}

Hence
\begin{equation*}
\begin{split}
(\star\star)&\ll\sup_{\substack{|u_1|, |u_2|,|u_3|\leq C^{\varepsilon}\\z\asymp \frac{\sqrt{X}}{C}}}\sum_{\substack{2\leq k\leq C^{\varepsilon}\mathcal{Z}\\ k\ \mathrm{even}}}\Gamma(k)\sum_{f\in\mathcal{B}_k}\bigg|\int_0^{\pi}e^{i(k-1)\xi}\sum_{b\asymp \mathcal{M}_1+\mathcal{M}_2}\sqrt{b}\,\overline{\rho_f(b)}\cdot e^{-i\sqrt{2b}z\sin\xi}\gamma^{\star}(b,z)\mathrm{d}\xi\bigg|^2\\
&\ll\int_0^{\pi}\sup_{\substack{|u_1|, |u_2|,|u_3|\leq C^{\varepsilon}\\z\asymp \frac{\sqrt{X}}{C}}}\sum_{\substack{2\leq k\leq C^{\varepsilon}\mathcal{Z}\\ k\ \mathrm{even}}}\Gamma(k)\sum_{f\in\mathcal{B}_k}\bigg|\sum_{b\asymp \mathcal{M}_1+\mathcal{M}_2}\sqrt{b}\,\overline{\rho_f(b)}\cdot e^{-i\sqrt{2b}z\sin\xi}\gamma^{\star}(b,z)\bigg|^2\mathrm{d}\xi\\
&\ll \mathcal{C}^{\varepsilon}\left(\mathcal{Z}^2+\mathcal{M}_1+\mathcal{M}_2\right)\sup_{z\asymp \frac{\sqrt{X}}{C}}\sum_{b}\big|\gamma^{\star}(b,z)\big|^2
\end{split}
\end{equation*}
by Lemma \ref{largesieve}.

In summary, we have arrived at
\begin{equation*}
\mathcal{HH}\ll\frac{C^{\varepsilon}HH'}{\Lambda}X^{-\frac{1}{2}}(\mathcal{M}_1\mathcal{M}_2)^{-\frac{1}{4}}\cdot \left(\mathcal{Z}^2+X\right)^{\frac{1}{2}}\|a\|_2\cdot \left(\mathcal{Z}^2+\mathcal{M}_1+\mathcal{M}_2\right)^{\frac{1}{2}}\sup_{z\asymp \frac{\sqrt{X}}{C}}\|\gamma^{\star}\|_2.
\end{equation*}

Now for our purpose it remains to deal with the $\ell^2$-norm $\|\gamma^{\star}\|_2=\left(\sum_{b}\big|\gamma^{\star}(b,z)\big|^2\right)^{\frac{1}{2}}$, where $\gamma^{\star}(b,z)$ is defined in \eqref{gammastar}. By Parseval, 

\begin{equation}
\begin{split}
\sum_{b\asymp \mathcal{M}_1+\mathcal{M}_2}\big|\gamma^{\star}(b,z)\big|^2\ll&
\int_{0}^{1}\left|\sum_{m_1\asymp\mathcal{M}_1}\lambda_{1}(m_1)\left(\frac{m_1}{\mathcal{M}_1}\right)^{-\frac{1}{4}+iu_1}e^{\pm iz\sqrt{m_1}}e(m_1\alpha)\right|^2\\
&\cdot\left|\sum_{m_2\asymp\mathcal{M}_2}\lambda_{2}(m_2)\left(\frac{m_2}{\mathcal{M}_2}\right)^{-\frac{1}{4}+iu_2}e^{\pm iz\sqrt{m_2}}e(m_2\alpha)\right|^2\mathrm{d}\alpha.
\end{split}
\end{equation}


\smallskip

Note that both $u_1$ and $u_2$ are of negligible size here. The sup-norm of the $m_2$-sum is bounded by $C^{\varepsilon}(\mathcal{M}_2^{1/2}+z\mathcal{M}_2)$, by Wilton's bound $\sum_{n\leq x}\lambda_f(n)e(\alpha n)\ll_f x^{1/2+\varepsilon}$ and partial summation. Next we open the square, getting
\begin{equation*}
\sum_{b}\big|\gamma^{\star}(b,z)\big|^2\ll C^{\varepsilon}(\mathcal{M}_2^{1/2}+z\mathcal{M}_2)^2\sum_{m_1\asymp\mathcal{M}_1}|\lambda_{1}(m_1)|^2\ll C^{\varepsilon}(\mathcal{M}_2^{1/2}+z\mathcal{M}_2)^2\mathcal{M}_1.
\end{equation*}

Similar sums without the twisted weight functions have appeared in the work of Jutila \cite{jut}. See also Blomer and Mili{\'c}evi{\'c} \cite[(8.5)]{bl-mi} for a similar sum.

\smallskip
Recall $\ \mathcal{M}_1\ll T_1=C^{\varepsilon}\frac{C^{2}X}{H^2}$, $\mathcal{M}_2\ll T_2=C^{\varepsilon}\frac{C^{2}X}{H'^2}$ $\mathcal{Z}=\frac{\sqrt{X(\mathcal{M}_1+\mathcal{M}_2)}}{C}\ll C^{\varepsilon}\frac{X}{H}$, $\Lambda\gg C^{2-\varepsilon}$, and $C=X^{1000}$ is a large parameter.
In particular, $\sup_{z\asymp \frac{\sqrt{X}}{C}}\sum_{b}\big|\gamma^{\star}(b,z)\big|^2\ll C^{\varepsilon}(\frac{X}{H'})^2\mathcal{M}_1\mathcal{M}_2$.
\smallskip
Then
\begin{equation*}
\begin{split}
\mathcal{HH}&\ll\frac{C^{\varepsilon}HH'}{C^2}X^{-\frac{1}{2}}(\mathcal{M}_1\mathcal{M}_2)^{-\frac{1}{4}}\cdot \left(\mathcal{Z}^2+X\right)^{\frac{1}{2}}\|a\|_2\cdot \left(\mathcal{Z}^2+\mathcal{M}_1+\mathcal{M}_2\right)^{\frac{1}{2}}\cdot \frac{X}{H'}\left(\mathcal{M}_1\mathcal{M}_2\right)^{\frac{1}{2}}\\
&\ll\frac{C^{\varepsilon}HH'}{C^2} X^{-\frac{1}{2}}\left(\frac{C^{2}X}{H^2}\cdot\frac{C^{2}X}{H'^2}\right)^{\frac{1}{4}}\left(\frac{X^2}{H^2}+X\right)^{\frac{1}{2}}\cdot \left(\frac{X^2}{H^2}+\frac{C^{2}X}{H^2}\right)^{\frac{1}{2}}\cdot\frac{X}{H'}\cdot \|a\|_2\\
&\ll C^{\varepsilon}(HH')^{\frac{1}{2}}\left(\frac{X^2}{H^2}+X\right)^{\frac{1}{2}}\cdot\frac{X^{1/2}}{H}\cdot\frac{X}{H'}\|a\|_2\\
&\ll C^{\varepsilon}\frac{X^{3/2}}{(HH')^{\frac{1}{2}}}\left(\frac{X}{H}+X^{\frac{1}{2}}\right)\|a\|_2.
\end{split}
\end{equation*}

By taking $H'=X/3$, we have $\mathcal{HH}\ll C^{\varepsilon}\frac{X}{H}\left(\frac{X}{H^{1/2}}+(XH)^{\frac{1}{2}}\right)\|a\|_2$, and hence
\begin{equation*}\sum_{X\leq n\leq 2X}a(n)E(n)\ll C^{\varepsilon}\frac{X}{H}\left(\frac{X}{H^{1/2}}+(XH)^{\frac{1}{2}}\right)\|a\|_2.
\end{equation*}

%



\subsection*{Acknowledgements.} The author is grateful to Valentin Blomer for several helpful comments and suggestions during the preparation of this article. He also thanks Sheng-Chi Liu for drawing his attention to \cite{blo} and Roman Holowinsky for his encouragement and comments.

\bigskip

%
%
%
%
%
%
%
%
%
%
%
%
%

\bibliographystyle{abbrv}
\bibliography{ref}

\end{document}